\documentclass[red,11pt,a4paper]{article}

\usepackage{cite}
\usepackage{comment}
\usepackage{amsmath,bm}
\usepackage{amscd}
\usepackage{amssymb}
\usepackage[pdftex]{color,graphicx,hyperref}
\usepackage[utf8]{inputenc}
\usepackage{latexsym}
\usepackage{color}
\usepackage{cases}
\usepackage{xfrac}
\usepackage[normalem]{ulem}
\usepackage{environ}
\usepackage{xcolor}
\usepackage{mathtools}
\usepackage{stackengine}

\usepackage{tikz}
\usetikzlibrary{positioning}

\bf
\usepackage{ika}
\renewcommand{\dd}{{\rm d}}

\newcommand{\T}{\mathbb{T}}

\newcommand{\vo}{{{\omega}}}

\newtheorem{defi}{Definition}[section]

\newtheorem{theo}[defi]{Theorem}
\newtheorem{cor}[defi]{Corollary}
\newtheorem{pro}[defi]{Proposition}

\newtheorem{rem}[defi]{Remark}

\NewEnviron{AS}{%
		\begin{align}
			\begin{split}
				\BODY
			\end{split}
		\end{align}
}
\title{Non-local dissipative Aw-Rascle model and its relation with Matrix-valued communication in Euler alignment}
\author{Nilasis Chaudhuri $^\dagger$
\and Jan Peszek $^\dagger$
\and Maja Szlenk $^{\dagger\ddagger}$
\and  Ewelina Zatorska \footnote{Corresponding author: ewelina.zatorska@warwick.ac.uk}
}

\date{\today}

\begin{document}
\maketitle

	\centerline{$^\dagger$ Institute of Applied Mathematics and Mechanics, University of Warsaw,}
	\centerline{ ul. Banacha 2, 02-097 Warszawa, Poland}
	\vspace{5mm}
\centerline{$^\ddagger$  Universit\'e Savoie Mont Blanc, CNRS, LAMA}

\centerline{73376 Le Bourget du lac, France}

\vspace{5mm}
\centerline{$^*$  Mathematics Institute, University of Warwick }

\centerline{Zeeman Building, Coventry CV4 7AL, United Kingdom}

\begin{abstract}
	
We compare the multi-dimensional generalisation of the Aw-Rascle model with the pressureless Euler-alignment system, in which the communication weight is matrix-valued. Our generalisation includes the velocity offset in the form of a gradient of a non-local density function, given by the convolution with a kernel $K$. We investigate connections between these models at the macroscopic, mesoscopic and macroscopic (hydrodynamic) level, and overview the results on the mean-field limit for various assumptions on $K$.

\end{abstract}

{\bf Keywords:}  Dissipative non-local Aw-Rascle system, Euler alignment system, matrix-valued communication weight, mean-field limit

{\bf MSC:}  35B40, 35L65, 35Q31

\bigbreak


\bigskip

\maketitle
\section{Introduction}
The purpose of this note is to explain the connection between the dissipative Aw-Rascle model and the Euler-alignment system with matrix-valued communication weight. We also review the well-posedness results for both of the models along with the literature concerning the derivation and the relations between their respective mico- meso- and macroscopic descriptions.  The main novelty of the paper is the adaptation of results from \cite{CC21} and \cite{PP23} to the dissipative non-local Aw-Rascle model, which we employ to investigate the micro-to-macroscopic mean-field limit, cf. Corollary \ref{corollary} and Theorem \ref{CCth}.\\
To begin with, let us recall the classical Aw-Rascle-Zhang  model (ARZ) of one lane road traffic, \cite{AR, Zhang}:
\[ \left\{ \begin{aligned}
		&\partial_t \vr + \partial_x(\vr u) = 0,  \\
		&\partial_t (\vr w) + \partial_x (\vr w u ) = 0,\\
  &{u} = {w} - P(\vr).
  \end{aligned} \right. \]
Similarly to compressible pressureless Euler equations, the ARZ model is a system of two conservation laws. The difference between the systems is in the  momentum equation, which for the ARZ model is associated with the the preferred velocity of motion $w$ and not the actual velocity $u$. The two velocities are related by the density-dependent offset function $P(\vr)$ which models formation of congestion on the road.

The multi-dimensional extension of this model, called the {\emph{dissipative Aw-Rascle model}},  was recently proposed in \cite{ABDM} to describe the flow of pedestrians, see also \cite{CGZ,CMPZ, CNPZ, M} for the relevant analytical results. It is again a system of conservation laws, where velocity vector fields $\vw$ and $\vu$ are related not by the scalar function $P(\vr)$, but by the gradient of the offset function $\Grad p(\vr)$.
In \cite{ABDM} a special singular function $p(\vr)$ was considered to take into account the volume exclusion effect and prevent overcrowding. However, the effects modelled by function $p(\vr)$ were of purely local character. 

In this note we  discuss the version of the dissipative Aw-Rascle model, in which the offset function is non-local and depends on the macroscopic density through convolution with the offset kernel $K(x)$
\begin{subnumcases}{\label{AR_multi}}
	\pt \vr+\Div (\vr \vu)=0,\label{AR_multi1}\\
	\pt(\vr \vw)+\Div(\vr \vu\otimes\vw)=0,\label{AR_multi2}\\
 \vu=\vw-\Grad K \ast  \vr.\label{def:w}
\end{subnumcases} 
We will refer to this system as the dissipative {\emph{non-local}} Aw-Rascle model, DNAR for short. We will review the mathematical results for the DNAR model and explain  the connection with Euler-alignment model with matrix-valued communication weight (EAM) derived below.

On the level of classical solutions, both systems: DNAR and EAM, are equivalent. Indeed, let us assume for now that system \eqref{AR_multi} has a classical solution in $\R^3$, and that the kernel $K$ is sufficiently smooth and symmetric $K(x)=K(-x)$.
Rewriting the momentum equation \eqref{AR_multi2} with indices, we obtain
\begin{align*}
\pt(\vr w_j)+\sum_{i=1}^3 \partial_{x_i} (\vr w_j u_i) =0, \; j=1,2, 3,
\end{align*}
with 
\[ w_j=u_j+ (\partial_{x_j}K ) \ast \vr . \] 
Then, using the continuity equation \eqref{AR_multi2}, we rewrite the momentum equation:
\begin{align*}
    0=&\pt(\vr u_j)+\sum_{i=1}^3 \partial_{x_i} (\vr u_j u_i) +	\pt(\vr (\partial_{x_j}K ) \ast \vr )+\sum_{i=1}^3 \partial_{x_i} (\vr (\partial_{x_j}K ) \ast \vr\,  u_i) \\
    &=\pt(\vr u_j)+\sum_{i=1}^3 \partial_{x_i} (\vr u_j u_i) + \vr \lr{ \partial_{x_j}K  \ast \pt \vr+\sum_{i=1}^3 (\partial_{x_i}\partial_{x_j}K ) \ast \vr\,  u_i  }.
\end{align*}
Again, with the help of the continuity equation and integration by parts, we deduce that
\begin{align*}
    \partial_{x_j}K  \ast \pt \vr & = -\int (\partial_{x_j}K)(x-y)  \lr{ \sum_{i=1}^3 \partial_{y_i} (\vr u_i) } \dy\\
    &=- \sum_{i=1}^3 \int (\partial_{x_i}\partial_{x_j}K ) (x-y) \vr(y) u_i(y) \dy .
\end{align*}
Consequently
\eqh{
    &\partial_{x_j}K  \ast \pt \vr+\sum_{i=1}^3 [(\partial_{x_i}\partial_{x_j}K ) \ast \vr]\,  u_i \\
    &=\sum_{i=1}^3 \lr{\int (\partial_{x_i}\partial_{x_j}K )(x-y) \lr{u_i(x)- u_i(y)} \vr(y) \dy. }
}
Therefore, coming back to the vector notation, we can rewrite system \eqref{AR_multi} in the form of the EAM system
\begin{subnumcases}{\label{EA}}
   	\pt \vr+\Div (\vr \vu)=0, \label{EA1} \\
\pt(\vr \vu )+ \Div(\vr \vu \otimes \vu)- \vr  \lr{\int \lr{\vu(y)- \vu(x)}^{\intercal}\boldsymbol{\Psi}(x-y)  \vr(y) \dy } =0, \label{EA2}
\end{subnumcases}
where 
\eq{\label{khess}
\boldsymbol{\Psi}=\nabla^2 K = \lr{\partial^2_{x_i x_j}K }_{i,j=1}^n.}
Notice that for the classical solution of \eqref{EA} and for a given sufficiently regular $\boldsymbol{\Psi}$, satisfying \eqref{khess}, the steps leading to derivation of system \eqref{EA} can be inverted. In this sense the EAM model \eqref{EA} is equivalent to the DNAR model \eqref{AR_multi}. Because both systems describe the evolution of macroscopic quantities (averaged density and velocity) characteristic for the models of fluids, we refer to them as {\emph{hydrodynamic}} or {\emph{macroscopic}} models. We summarise the above argument in the form of the following proposition.

\begin{pro}\label{equiv}
    Assume that $K$ is sufficiently smooth and symmetric i.e. $K(x)= K(-x)$. Then classical solutions $(\vr,\vw)$ to the DNAR system \eqref{AR_multi} are equivalent, through the transformation \eqref{def:w}, to classical solutions to the EAM system \eqref{EA}.
\end{pro}

\begin{figure*}[t]\label{fig}
    \centering
    
    \scriptsize{
\begin{tikzpicture}[
  node distance = 2cm,
  box/.style = {draw, minimum height=7mm, align=center,
  }]
  \node (n1)   {microscopic};
  \node[right=0.5cm of n1] (n2)  [box]   {$\begin{aligned}
          \dot{x}_i &= v_i, \\
          \dot{\omega}_i &= 0, \\
          v_i &= \omega_i-\frac{1}{N}\sum_{i=1}^N\nabla K(x_i-x_j)
      \end{aligned}
  $};
  \node[right=3.5cm of n2] (n3)  [box]   {$\begin{aligned}
      \dot{x}_i &= v_i, \\
      \dot{v}_i &= -\frac{1}{N}\sum_{i=1}^N \nabla^2K(x_i-x_j)(v_i-v_j)
  \end{aligned}$};

  \node[below=4cm of n1] (n4)  {mesoscopic};
  \node[right=0.5cm of n4] (n5) [box] {$\begin{aligned}
      &\partial_t\mu + \mathrm{div}_x({u}[\varrho_t]\mu) = 0, \\
      &{u}[\varrho_t](x,{\omega}) = {\omega}-\nabla K\ast\varrho_t, \\
      &\varrho_t(x) = \int_{\mathbb{R}^d}\;\mu_t(x,\mathrm{d}\omega)
  \end{aligned}$ };
  \node[right=3.5cm of n5] (n6) [box] {$\begin{aligned}
      &\partial_t\mu + {v}\cdot\nabla_x\mu - \mathrm{div}_{{v}}(F(\mu)\mu) = 0, \\
      &F(\mu)(t,x,{v}) = \int_{\mathbb{R}^d}\nabla^2K(x-y)({v}-\Tilde{{v}})\;\mu(t,y,\mathrm{d}\tilde{{v}})\;\mathrm{d}y
  \end{aligned}$};
  \node[below=4cm of n4] (n7)  {macroscopic};
  \node[right=0.5cm of n7] (n8) [box] {$\begin{aligned}
      &\partial_t\varrho + \mathrm{div}(\varrho\boldsymbol{u}) = 0, \\
      &\partial_t(\varrho\boldsymbol{w}) + \mathrm{div}(\varrho\boldsymbol{u}\otimes\boldsymbol{w}) = 0, \\
      &\boldsymbol{u} =\boldsymbol{w}-\nabla K\ast\varrho
  \end{aligned}$};
  \node[right=2.5cm of n8] (n9) [box] {$\begin{aligned}
      &\partial_t\varrho + \mathrm{div}(\varrho\boldsymbol{u}) = 0, \\
      &\partial_t(\varrho\boldsymbol{u}) + \mathrm{div}(\varrho\boldsymbol{u}\otimes\boldsymbol{u}) = -\varrho\int_{\mathbb{R}^d}\nabla^2K(x-y)(\boldsymbol{u}(x)-\boldsymbol{u}(y))\varrho(y)\;\mathrm{d}y
  \end{aligned}$};
  \draw[stealth-stealth,shorten >=4pt,shorten <=4pt] (n2) -- (n3) node[pos=0.5,above]{$K$ twice differentiable};
  
  \draw[-latex,shorten >=4pt,shorten <=4pt] (n2.south)
     -- (n2.south |- n5.north) 
     node[pos=0.5,right]{Theorem \ref{whelp}};
  \draw[-latex,shorten >=4pt,shorten <=4pt] (n3.east)
     -| (n9.east |- n9.north) 
     node[pos=0.6,left]{Theorem \ref{CCth}};
  \draw[latex-,shorten >=4pt,shorten <=4pt] (n5.south)
     -- (n5.south |- n8.north) 
     node[midway,right]{Proposition \ref{prop:monokin}};


  \draw[stealth-stealth,shorten >=4pt,shorten <=4pt] (n8) -- (n9) node[pos=0.5,above]{Proposition \ref{equiv}};

  \draw[stealth-stealth,shorten >=4pt,shorten <=4pt, color=red, dashed] (n5) -- (n6) node[pos=0.5, above,color=black]{\cite{Kim,PP22}};
\end{tikzpicture}
}

    \caption{Microscopic, mesoscopic and macroscopic formulations of the DNAR model (left) and of the EAM model (right) and relations between the systems discussed in this article. With a dashed red arrow we indicate the relations described in other papers.}
    \label{Fig1}
\end{figure*}
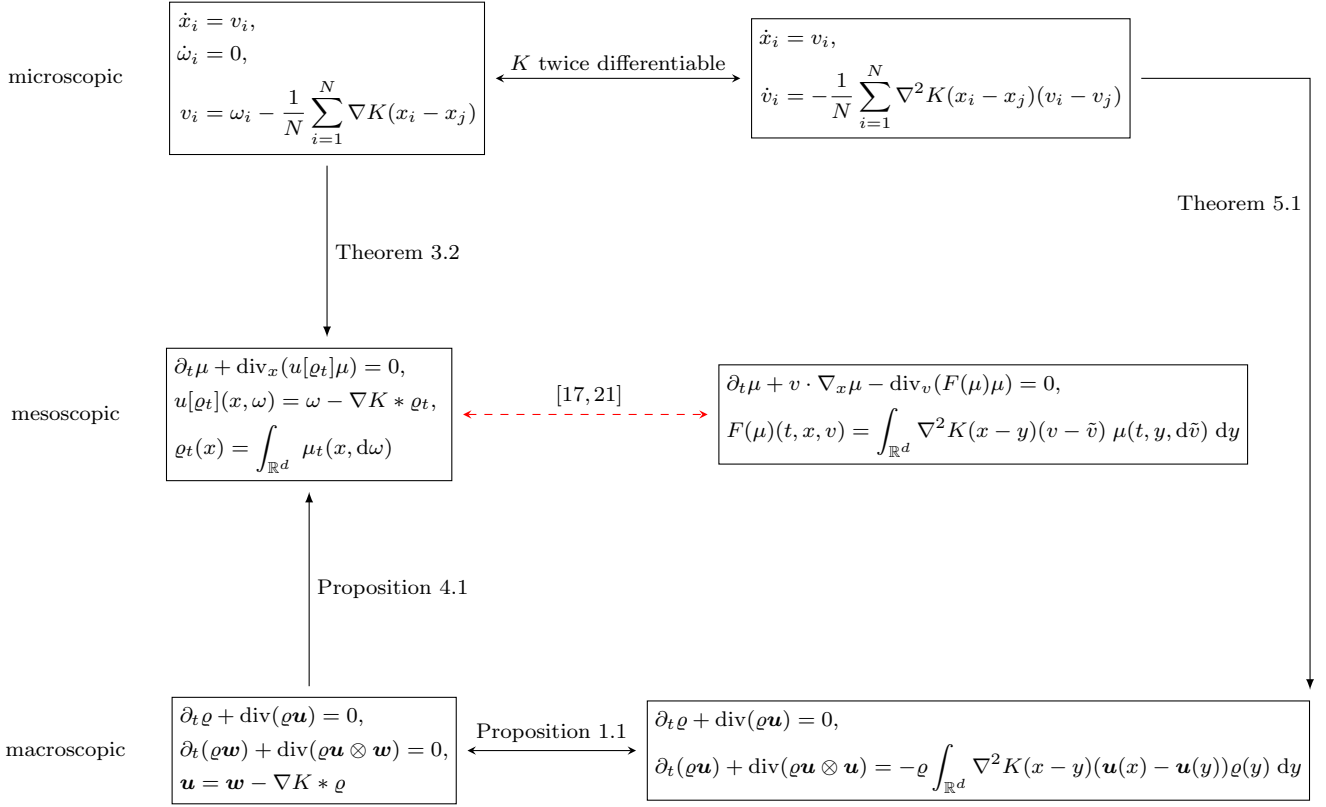

    Taking as a point of reference Figure \ref{fig}, we note that the microscopic equivalence between the DNAR and EAM systems follows by a direct differentiation in time of the former. In one space dimension, it has been the basis for cluster predictability for the Cucker-Smale model,  see \cite{HPZ, HKPZ}. In the case of regular kernels $K$, the equivalence has been extended by Kim \cite{Kim21} to the kinetic and hydrodynamic levels in one space dimension. Choi and Zhang \cite{CZ} used it as a means to establish the existence of solutions to the Cucker-Smale kinetic equation for weakly singular communication weight. In alignment dynamics, weakly singular communication weight is associated with the kernel $K(x)=|x|^{2-\alpha}$ with $\alpha\in(0,1)$, which translates to the singular communication weight $D^2K(x)\approx |x|^{-\alpha}$, on the level of the Cucker-Smale model. The equivalence relation between higher-dimensional Cucker-Smale-type models and their first-order reduction based on the DNAR dynamics was then studied by Kim \cite{Kim} and extended to the kinetic models for smooth and strictly convex potentials $K$. In the case of singular communication (again, $K(x)=|x|^{2-\alpha}$) the equivalence was obtained by the second author and Poyato in \cite{PP22}.
    The derivation of the hydrodynamic DNAR model and the EAM model, as well as their equivalence, were never studied as far as we are aware. Also, the case of more general interaction potentials $K$ is widely open. It should be noted however that the hydrodynamic Euler-alignment system (with classical scalar communication) has been derived rigorously, cf. \cite{FK, KMT, FP, HT}.

In this note, we focus only on the multi-dimensional models, because in 
1D the matrix-valued communication weight becomes the usual scalar one, and the system \eqref{EA} reduces to relatively  well-studied Euler-alignment system, see for instance  \cite{kis, tad1,tad2, tad4, LS}.
Most studies of the 1D Euler-alignment system rely to some extent on the quantity
\eqh{
e(t,x):= 
 \px u(t,x) + \int_{\Omega}\psi(x-y)(\vr(t,x)-\rho(t,y))\dy,\quad \psi(x)=\px^2 K(x)=\frac{1}{|x|^\alpha},
}
where $\Omega$ is either the whole line $\R$ or the torus $\T$.
The key property of $e$ is that it satisifes the continuity equation
\eqh{
\pt e + \px(ue) = 0.
}
From the DNAR perspective, noting that $\px w=\px u+\px^2 K\ast \vr$, it translates into the continuity equation
\eqh{
\pt  (\px w)+\px(u\px w)=0,
}
at least on the support of $\vr$, which is a natural consequence of the transport equation for $w$. In other words, the $e$ quantity is nothing but the spacial derivative of the preferred velocity $w$.
The algorithm of obtaining higher order entropies was discussed in \cite{CDS20}, and it was used as a tool to propagate initial regularity of classical solutions (provided no vacuum has formed) instead of the classical energy method.
Equation \eqref{AR_multi2} can be therefore seen as a natural generalization of this concept to the multi-dimensional setting, for the entropy of order ``-1''. Indeed, the (formally) derived equation for  $\Div\vw$,
\eqh{\pt  (\Div \vw)+\Div(\vu\cdot\Grad \vw)=0,}
no longer has the structure of the continuity equation.
It is therefore not clear how to generate a strong solution to any of the multidimensional systems \eqref{EA} or \eqref{AR_multi}.

Throughout  our survey we will investigate the equivalence between the systems \eqref{AR_multi} and \eqref{EA}, as well and their {\emph{microscopic}} and {\emph{mesoscopic}} (also called {\emph{kinetic}}) counterparts. The corresponding mathematical models, as well as selected relations between them, are depicted on the schematic diagram in Fig. \ref{Fig1}.

The paper is organized as follows.  First, in Section \ref{sec:micromeso}, we perform the mean-field limit for DNAR system, under the assumption that $K$ is $\lambda$-convex or weakly singular (as in (\ref{defK})). Theorem \ref{whelp} provides the existence and uniqueness result to mesoscopic formulation of the DNAR system in a certain class of weak measure-valued solutions, and shows that the solution is a mean-field limit of the microscopic ODE system. The connection between the macroscopic (hydrodynamic) DNAR system (\ref{AR_multi}) and its kinetic counterpart is explained in Section \ref{sec:mesomacro}, in Proposition \ref{prop:monokin}. It shows that a suitably regular solution to (\ref{AR_multi}) generates a monokinetic solution to the corresponding mesoscopic system. 
In Corollary \ref{corollary} we synthesise the results from the two previous sections and provide (conditional) results for the convergence of the empirical measures to the weak solutions of the macroscopic DNAR system.
 Lastly, in Section \ref{micromacro}, we present the direct micro-to-macroscopic mean-field limit for strong solution of the macroscopic system (\ref{EA}), see Theorem \ref{CCth}. This is a matrix-valued generalization of the result from \cite{CC21}, that requires stronger assumptions on the potential, namely  $\Psi=\nabla^2K\in W^{2,\infty}(\R^d)$. The existence of (local-in-time) classsical solutions requires even stronger assumption,  $\boldsymbol{\Psi}\in C^1_c(\R^d)$, as stated in Theorem \ref{Th:local}.

\section{Notation}

Throughout the paper by ${\mathcal P}(\R^d)$ we denote the space of probability measures on $\R^d$ and by ${\mathcal P}_p(\R^d)$, with $p\geq 1$ -- the space of probability measures with finite $p$-moment. We say that $\mu_n\to\mu$ narrowly, if and only if
$$ \int_{\R^d}\phi \,{\rm d} \mu_n \to \int_{\R^d}\phi \,{\rm d}\mu$$
for all $\phi\in C_b(\R^d)$ (i.e. all bounded-continuous functions $\phi$). By ${\mathcal P}(\R^d)-narrow$ we denote the space of probability measures with the narrow topology; we specifically use the space $C([0,T];{\mathcal P}(\R^{2d})-narrow)$ as the space of all continuous mappings from $[0,T]$ to the ${\mathcal P}(\R^d)-narrow$ space.

 For a Borel-measurable mapping $F:\R^d \rightarrow \R^n$ and a probability measure $\mu$ we define the \textit{pushforward measure} $F_{\#}\mu$ by
    \begin{equation*}
    F_{\#}(\mu)(B):= \mu\left( F^{-1}(B) \right) \quad \text{for Borel sets } B.
    \end{equation*}
    The main property of the pushforward measure,  used throughout the paper, is the change of variables formula
    \begin{equation*}
        \int_{\R^n} g \,{\rm d} (F_{\#}\mu) = \int_{\R^d} g \circ F \,{\rm d} \mu,
    \end{equation*}
    whenever $g\circ F$ is $\mu$-integrable.

We shall employ distances based on optimal transport commonly associated with the ${\mathcal P}_p(\R^d)$ spaces. These include the bounded-Lipschitz distance
\eqh{ 
d_{BL}(\mu,\nu) = \sup\left\{\left|\int_{\R^d}\phi \,{\rm d}\mu - \int_{\R^d}\phi \,{\rm d}\nu\right|:\quad \phi \mbox{ is Lipschitz continuous with } [\phi]_{Lip}\leq 1, |\phi|_\infty\leq 1\right\}, }
where by $[\phi]_{Lip}$ we denote the Lipschitz constant of $\phi$. Another metric used throughout the paper is the Wasserstein distance
$$ W_2(\mu,\nu) = \inf_{\gamma\in \Gamma(\mu,\nu)} \int_{\R^d\times\R^d}|x-y|^2\,{\rm d}\gamma(x,y),$$
where the infimum is taken over the set $\Gamma(\mu,\nu)$ of admissible plans, i.e. over all probability measures $\gamma\in {\mathcal P}(\R^d\times\R^d)$ with $x$- and $y$-marginals equal to $\mu$ and $\nu$, respectively. In other words 
$$ \Gamma(\mu,\nu) = \{\gamma\in {\mathcal P}(\R^d\times\R^d):\quad (\pi_x)_\#\gamma=\mu \mbox{ and } (\pi_y)_\#\gamma=\nu\}, $$
where $\pi_x$ and $\pi_y$ are projections to the first and the second component, respectively.

Throughout the paper we write
$$ \int_{\R^{d}} f(x,y) \mu(x,\dd y) $$
to indicate that the integration above is taken only with respect to $y$. Integrating with respect to all available variables we use the simplified notation
$$ \int_{\R^{2d}} f(x,y)\dd\mu .$$

\section{From micro- to mesoscopic DNAR dynamics}\label{sec:micromeso}

In full analogy to the derivation of AR system from the Follow-the-Leader particle dynamics, one may try to connect the macroscopic equations \eqref{AR_multi} with their microscopic counterpart. 
For the dissipative AR system with offset function $p(\vr)= (\vr^{-1}-\vr_{{\rm{max}}}^{-1})^{-\gamma}$ the formal procedure was explained in \cite{ABDM}. In this section we discus various rigorous results concerning derivation of the kinetic variant of DNAR model from the system of interacting active particles, for which the offset function is non-local. 



The microscopic analogue of the DNAR system, describing interactions of $N$ agents with positions $x_i(t)\in\R^d$ and velocities $ v_i(t)\in\R^d$ at time $t\geq0$, is given by the ODE system:
\begin{subnumcases}{\label{ar}}
\dot{x}_i=v_i, \label{ar1}\\
\dot{\omega}_i=0,\label{ar1a}\\
v_i=\omega_i-\displaystyle\frac{1}{N}	\displaystyle\sum_{j=1}^N\Grad K(x_i-x_j).\label{ar2}
\end{subnumcases}
Here, $v_i$ is the actual velocity of the $i$th particle, while $\vo_i$ is the desired velocity of motion.
Note that, differentiating  equation \eqref{ar2} w.r.t. $t$, we can rewrite the system in the form of the {\emph{Cucker-Smale}-type} system of alignment dynamics
\begin{subnumcases}{\label{cs}}
\dot{x}_i=v_i,\label{cs1}\\
\dot{v}_i=\displaystyle\frac{1}{N}	\displaystyle\sum_{j=1}^N(v_j-v_i)\boldsymbol{\Psi}(x_i-x_j),\label{cs2}
\end{subnumcases}
where $\boldsymbol{\Psi}$ is related to $K$ by \eqref{khess}. More specifically, \eqref{cs} is a variant of the Cucker-Smale model with matrix-valued communication as in \cite{ST21}.  In full analogy with the equivalence between the hydrodynamic models, cf. Proposition \ref{equiv}, their microscopic counterparts, systems \eqref{ar} and \eqref{cs}, are also equivalent as long as the existence classical solutions is guaranteed by the Cauchy–Lipschitz theory.

In the limit $N\to\infty$, one formally describes system \eqref{ar} by the probability $\mu_t(x,\vo)$ of finding agents with desired velocity $\vo$ around point $x$ at time $t$. It is subject to the kinetic equation
\begin{subnumcases}{\label{ark}}
\pt \mu+\Div_x(u[\vr_t]\, \mu)=0, \label{ark1}\\
 u[\vr_t](x,\omega)=\omega-\Grad_x K\ast \vr_t(x), \label{ark1a}\\
 \vr_t(x)=\int_{\R^d} \mu_t(x,{\rm d}\omega) \label{ark2}.
\end{subnumcases}

It was shown in \cite{PP23} that an appropriate weak measure-valued solution to \eqref{ark} is a gradient flow of the interaction energy
\eq{\label{defW}
{\mathcal W}(\mu) = \int_{\R^{4d}}K(x-x'){\rm d}[\mu\otimes\mu],
}

\noindent
under one of the following assumptions
\begin{equation}\label{defK}
    \begin{split}
        K(x)&=K(-x) \mbox{ is lower semicontinuous and } \lambda-\mbox{convex for } \lambda>0, \quad \mbox{or}\\
        K(x)&=\frac{1}{2-\alpha}\frac{1}{1-\alpha} |x|^{2-\alpha}, \quad x\in \R^d,\quad \alpha\in(0,1).
    \end{split}
\end{equation}
In \cite{PP23} kernel $K$ in \eqref{defK} is refereed to as ``weakly singular'' due to the singularity of the associated alignment communication kernel $\boldsymbol{\Psi}\approx|x|^{-\alpha}$ in \eqref{cs2}.

 Using the gradient flow interpretation of \eqref{ark} leads to a variety of results such as well-posedness, uniform-in-time contractivity. It is also used to rigorously derive system \eqref{ark} from the underlying particle dynamics governed by \eqref{ar}. 
The definition of the solution  and the statement of some main results from \cite{PP23} are recalled below.

\begin{defi}[Weak formulation for \eqref{ark}]\label{Def:PP}
Consider any $T>0$ and $\mu_0\in{\mathcal P}_2({\mathbb{R}^{2d})}$. We say that  $\boldsymbol{\mu}\in C([0,T],\mathcal{P}(\mathbb{R}^{2d})-\mbox{narrow})$, satisfying the assumption
\begin{equation}\label{E-weakeqF}
\int_0^T\int_{\mathbb{R}^{2d}}|x|^{1-\alpha}\,{\rm d}\mu_t(x,\omega)\,{\rm d}t<\infty,
\end{equation}
is a \emph{weak measure-valued solution} to \eqref{ark} in the time interval $[0,T]$, subject to initial datum $\mu_0$, if the following weak formulation is verified
\begin{align}\label{weakeqF}
- \int_{\mathbb{R}^{2d}} \eta_0\,{\rm{d}}\mu_0 = \int_0^T\int_{\mathbb{R}^{2d}}\Big(\partial_t\eta_t(x,\omega) + (\omega-\nabla K*\vr_t(x))\cdot\nabla_x\eta_t(x,\omega)\Big) \,{\rm d}\mu_t \,{\rm d}t,
\end{align}
for all test functions $\eta \in C^1_b([0,T]\times\R^{2d})$, compactly supported in $[0,T)$.
\end{defi}
In \cite{PP23} system \eqref{ark} has been reformulated as a gradient flow with respect to a fibered Wasserstein-2 distance. Indeed using the disintegration theorem, cf. Theorem \ref{T:disint} in the appendix, each time-slice $\mu_t$ of any solution $\mu$ to \eqref{ark} can be disintegrated with respect to its $\omega$-marginal, called $\nu$, and represented as
\begin{equation}\label{disint}
    \int_{\R^{2d}}\psi(x,\omega)\,{\rm d}\mu_t(x,\omega) = \int_{\R^{d}} \left(\int_{\R^d}\psi(x,\omega)\,{\rm d}\mu_t^\omega(x)\right)\,{\rm d}\nu(\omega),
\end{equation}
for all Borel-measurable functions $\psi:\R^{2d}\longrightarrow [0,+\infty)$. The key observation is that the $\omega$-marginal $\nu$ does not change in time and system \eqref{ark} can be disintegrated leading to a fiber-by-fiber family of continuity equations
\begin{align*}
&\partial_t\mu^\omega+{\rm div}_x({u}[\vr_t](\cdot,\omega)\,\mu^\omega)=0,\quad t\geq 0,\ x\in \mathbb{R}^d,\\
&{u}[\vr_t](x,\omega)=\omega-\int_{\mathbb{R}^d}\nabla K*\mu_t^{\omega'}(x)\,{\rm d}\nu(\omega'),\\
&\mu_{t=0}^\omega=\mu_0^\omega,
\end{align*}
for $\nu$-a.e. $\omega\in \mathbb{R}^d$. We remark that the interplay between different $\omega$-fibers of $\{\mu_t^\omega\}_{\omega\in\R^d}$ is highly complex since the velocity field ${u}[\vr_t]$ is reconstructed as an average of the combined effects of all fibers simultaneously. The natural space for this problem, used in \cite{PP23}, is the fibered Wassertein space
\begin{equation}\label{dnu}
\begin{split}
\mathcal{P}_{2,\nu}(\mathbb{R}^{2d})&:=\left\{\mu\in \mathcal{P}(\mathbb{R}^{2d}),\ (\pi_{\omega})_\#\mu=\nu:\,\int_{\mathbb{R}^{2d}}\vert x\vert^2\,{\rm d}\mu(x,\omega)<\infty\right\},\\
W_{2,\nu}(\mu_1,\mu_2)&:=\left(\int_{\mathbb{R}^d}W_2^2(\mu_1^\omega,\mu_2^\omega)\,{\rm d}\nu(\omega)\right)^{1/2},
\end{split}
\end{equation}

\noindent
in which $\mu$ becomes a gradient flow of the interaction energy ${\mathcal W}$ in \eqref{defW}.
Such an interpretation leads to the following results.
\begin{theo}[Well-posedness and the mean-field limit]\label{whelp}
    Suppose that $K$ satisfies \eqref{defK}. Let $\nu\in {\mathcal P}(\R^d)$ be the given distribution of $\omega$. For any initial datum $\mu_0\in {\mathcal P}_{2,\nu}(\R^{2d})$ there exists a unique solution to \eqref{ark} in the sense of Definition \ref{Def:PP}.

Moreover the solution $\mu$ is a mean-field limit, in $W_2$ metric, of atomic solutions
\begin{equation}\label{E-empirical-measures}
\mu_t^N:=\frac{1}{N}\sum_{i=1}^N\delta_{x_i^N(t)}(x)\otimes \delta_{\omega_i^N}(\omega),
\end{equation}
where $(x_i,\omega_i)_{i=1}^N$ satisfy system \eqref{ar}.
\end{theo}
\begin{theo}[Uniform contractivity in $W_{2,\nu}$]\label{T-contractivity-W2nu}
Suppose that $K$ satisfies \eqref{defK}. Consider any $\nu\in \mathcal{P}(\mathbb{R}^d)$ and any compactly supported initial data $\mu_0^1,\mu_0^2\in \mathcal{P}_{2,\nu}(\mathbb{R}^{2d})$ with the same center of mass, {\it i.e.},
\begin{equation*}
\int_{\mathbb{R}^{2d}}x\,{\rm d}\mu_0^1(x,\omega)=\int_{\mathbb{R}^{2d}} x\,{\rm d}\mu_0^2(x,\omega).
\end{equation*}
Let ${\mu}^1$ and ${\mu}^2$ be the weak measure-valued solutions to \eqref{ark} issued at $\mu_0^1$ and $\mu_0^2$. Then,
$$W_{2,\nu}(\mu_t^1,\mu_t^2)\leq e^{-2c_0t} \,W_{2,\nu}(\mu_0^1,\mu_0^2),$$
for each $t\geq 0$. Here, either $c_0=\lambda$ if $K$ is $\lambda$-convex, or $c_0 = D_0^{-\alpha}$,
where $D_0\in \mathbb{R}_+$ is a positive constant related to a uniform bound on the spacial spread of the supports of $\mu^1$ and $\mu^2$.
\end{theo}

\begin{cor}[Convergence to equilibrium]
Consider any $\nu\in \mathcal{P}(\mathbb{R}^d)$ verifying 
$$\omega_c:=\int_{\mathbb{R}^d}\omega\,{\rm d}\nu(\omega)=\int_{\mathbb{R}^{2d}}\omega\,{\rm d}\mu_0(x,\omega)=0,$$
any compactly supported initial datum $\mu_0\in \mathcal{P}_{2,\nu}(\mathbb{R}^{2d})$ and let ${\mu}$ be the weak measure-valued solution of \eqref{ark} issued at $\mu_0$. Then, there exists a unique compactly supported equilibrium $\mu_\infty\in \mathcal{P}_{2,\nu}(\mathbb{R}^{2d})$ of \eqref{ark}, such that
\begin{equation*}
\int_{\mathbb{R}^{2d}}x\,{\rm d}\mu_0(x,\omega)=\int_{\mathbb{R}^{2d}}x\,{\rm d}\mu_\infty(x,\omega).
\end{equation*}
In addition, we obtain
$$W_{2,\nu}(\mu_t,\mu_\infty)\leq e^{-2c_0t} \,W_{2,\nu}(\mu_0,\mu_\infty),$$
for any $t\geq 0$, where $c_0$ is as in Theorem \ref{T-contractivity-W2nu}.
\end{cor}

Let us discuss the above results. First, observe that the fibered Wasserstein distance $W_{2,\nu}$, being stronger than the classical $W_2$ metric, i.e. $ W_2(\mu_1,\mu_2)\leq C W_{2,\nu}(\mu_1,\mu_2)$, is tailored to work well with heterogeneus continuity equation such as \eqref{ark}. It takes into the account that, following \eqref{ar1a}, the motion occurs only along the $\omega$-fibers of $\mu$, and assigns an infinite cost of moving mass between  $\omega$-fibers.
On the other hand it is completely unsuited for the purpose of the mean-field limit, since the $\omega$-marginals of $\mu^N$ do not agree with one another or with the $\omega$-marginal of $\mu$. To put it simply, formula \eqref{dnu} does not make sense. Another metric is needed, which enables comparison of measures with different $\omega$-marginals. In \cite{PP23} the authors use the {\it adapted} (also known as {\it nested}) Wasserstein distance defined for any $\mu_1,\mu_2\in \mathcal{P}_2(\mathbb{R}^{2d})$, by the formula
\begin{equation}\label{E-adapted-Wassersein-distance}
AW_2(\mu_1,\mu_2):=\left(\inf_{\hat{\nu}\in \Gamma(\nu_1,\nu_2)}\int_{\mathbb{R}^{2d}} \left(W_2^2(\mu_1^\omega,\mu_2^{\omega'})+\vert \omega-\omega'\vert^2\right)\,{\rm d}\hat{\nu}(\omega,\omega')\right)^{1/2},
\end{equation}
where $\nu_1:=\pi_{\omega\#}\mu_1$ and $\nu_2:=\pi_{\omega\#}\mu_2$, and $\{\mu_1^\omega\}_{\omega\in \mathbb{R}^d}$ and $\{\mu_2^{\omega'}\}_{\omega'\in \mathbb{R}^d}$ are the associated families of disintegrations established in \eqref{disint}.

In \cite{PP23} a contractivity result akin to Theorem \ref{T-contractivity-W2nu} but with $AW_2$ was shown. Then it was used to derive the mean-field limit in the $AW_2$ metric. For simplicity\footnote{For example, it is unclear if the approximation of general probability measures with empirical measures is possible in the $AW_2$ metric.}, we formulate the results only in the $W_2$ and $W_{2,\nu}$ distances, the finer estimates in $AW_2$ can be found in the original work. Of course $AW_2$ is stronger than the classical $W_2$ Wasserstein distance i.e.
    $ W_2(\mu_1,\mu_2)\leq C AW_2(\mu_1,\mu_2).$
    Moreover if the $\omega$-marginals of $\mu_1$ and  $\mu_2$ agree and are equal to $\nu$, we have
    $ AW_2(\mu_1,\mu_2) = W_{2,\nu}(\mu_1,\mu_2), $
    which we infer by taking $\hat{\nu} = \delta_{\omega}\otimes\nu(\omega)$ in \eqref{E-adapted-Wassersein-distance}.

\section{From meso- to macroscopic DNAR dynamics}\label{sec:mesomacro}
In order to derive the DNAR system \eqref{AR_multi} from the corresponding kinetic equations \eqref{ark} we  proceed as in the case of derivation of compressible Euler equations from Vlasov-type equations, i.e. by taking the $\omega$-moments of equation \eqref{ark1}.

Let us define the macroscopic variables
\eqh{\vr=\int_{\R^d}  \mu(x, {\rm d}\vo),\quad \vr\vw=\int_{\R^d}\vo \;\mu(x, {\rm d}\vo).
}
We first integrate equation \eqref{ark1a}  w.r.t. $\vo$ to deduce that
\eqh{
\int_{\R^d} u[\vr_t]\, \mu(x, {\rm d}\vo)=\vr\vw-\vr\Grad K\ast\vr.}
Next, integrating equation \eqref{ark1} we obtain the continuity equation in the form:
\eqh{
\pt\vr+\Div(\vr\vw-\vr\Grad K\ast\vr)=0,
}
and so, defining as in \eqref{def:w} 
\eq{\label{new_defu}
\vu:=\vw-\Grad K\ast\vr,}  
we derive the continuity equation  \eqref{AR_multi1} of the DNAR system.

Next, we compute the first moment of \eqref{ark1} in $\vo$, obtaining
\eq{\label{intw}
\pt(\vr\vw)+ \Div\lr{\int_{\R^d} u[\vr_t]\otimes\vo \, \mu(x, {\rm d} \vo)}=0.
}

\noindent
Equivalently, using \eqref{ark1a} and \eqref{new_defu} subsequently, we can write
\eqh{
\int_{\R^d}u[\vr_t]\otimes\vo \, \mu(x, {\rm d} \vo)&=
\int_{\R^d}\lr{\omega-\Grad_x K\ast \vr_t(x)}\otimes\vo \, \mu(x, {\rm d} \vo)\\
&=\int_{\R^d}\,\omega\otimes\vo \, \mu(x, {\rm d} \vo)-\vr_t(\Grad_x K\ast \vr_t(x))\otimes\vw\\
&=\int_{\R^d}\,\omega\otimes\vo \, \mu(x, {\rm d} \vo)+\vr\vu\otimes\vw-\vr\vw\otimes\vw.
}
As a consequence, we can rewrite \eqref{intw} in the form
\eqh{
\pt(\vr\vw)+ \Div(\vr\vu\otimes\vw)+\Div \mathbf{\Pi}=0,
}
where
\eq{\label{pres}\mathbf{\Pi}=\int_{\R^d} (\vo-\vw)\otimes \vo\, \mu(x,{\rm d} \vo) .}

In this sense, the DNAR system is similar to the compressible Euler equation, with transported velocity $\vu$ replaced by $\vw$ and the monokinetic ansatz, characterising the pressureless case, translated to
\eqh{\mu(x,\vo)=\vr(t,x) \delta(\vo-\vw).}
Rigorous justification of hydrodynamic limit of the kinetic Cucker-Smale model to the pressureless Euler system with nonlocal alignment force was obtained by Figalli and Kang \cite{FK}. It was shown recently in \cite{FP} that in the strongly singular case, every weakly continuous and compactly supported solution to the kinetic Cucker-Smale equation is monokinetic, and its macroscopic variables satisfy the Euler-alignment system without any pressure term akin to  \eqref{pres} and with no necessity to take a hydrodynamic limit as in \cite{FK}.

The rigorous justification of this closure for general $K$ will be the purpose of our study in the future. In this note we restrict ourselves  to showing the following result  for weakly singular case. 

\begin{pro}\label{prop:monokin}
Suppose that $K$ is as in \eqref{defK} and in the $\lambda$-convex case it is $C^1$. Let $0\leq\vr\in C([0,T];{\mathcal P}(\R^d)-\mbox{narrow})$ and  $\vw\in C^1_b([0,T]\times\R^d)$ be a weak solution to \eqref{AR_multi} (tested with $C^1_b$ functions). If
\begin{equation}\label{eq:p2dlarw}
     \int_{\R^d}|x|^2 + |\vw_0(x)|^2 \,{\rm d}\vr_0< +\infty 
\end{equation}
and
\begin{equation}\label{eq:a-1dlarw}
     \int_0^T\int_{\R^d}|x|^{1-\alpha}{\rm d}\vr_t(x) \, {\rm d} t <+\infty,
\end{equation}
then the measure
\eqh{\mu_t(x,\vo)=\vr_t(x) \delta(\vo-\vw)}
is a measure solution to \eqref{ark} in the sense of Definition \ref{Def:PP}.


\end{pro}

\pf Let us check all the assumptions required in Definition \ref{Def:PP}. The fact that $\mu_0\in{\mathcal P}_2(\R^{2d})$ follows from \eqref{eq:p2dlarw}. To check the narrow continuity of $\mu$, let $\psi$ be a bounded-Lipschitz function\footnote{Narrow convergence is defined as tested by bounded-continuous functions, but by Portmanteau theorem it is equivalent to testing with bounded-Lipschitz functions.} and suppose that $t_n\to t$. We have
    \begin{align*}
      &  \left|\int_{\R^{2d}}\psi(x,\vo)\,{\rm d} \vr_{t_n}\delta(\vo-\vw_{t_n}(x)) - \int_{\R^{2d}}\psi(x,\vo) \,{\rm d} \vr_{t_n}\delta(\vo-\vw_t(x))\right|\\
       & \leq   \int_{\R^{d}}|\psi(x,\vw_{t_n}(x)) - \psi(x,\vw_{t}(x))| \,{\rm d} \vr_{t_n}(x) + \left|\int_{\R^{d}}\psi(x,\vw_t(x)) \,{\rm d} (\vr_{t_n}(x) -\vr_{t}(x))\right| =: I_n + II_n.
    \end{align*}
    By the Lipschitz continuity of $\psi$ (with constant $L_\psi$) as well as the regularity assumption on $\vw$, we have
    $$ I_n\leq L_\psi|\vw_{t_n}-\vw_t|_\infty\xrightarrow{n\to\infty} 0, $$
    and by the narrow continuity of $\vr$ we also have $II_n\to 0$, since $\psi(x,\vw_t(x))$ is a bounded-continuous function.

    Since condition \eqref{E-weakeqF} follows immediately from \eqref{eq:a-1dlarw}, 
    it remains to show the weak formulation \eqref{weakeqF}. To such an end let $\eta\in C^1_b([0,T]\times\R^{2d})$ and let $\widetilde{\eta}(t,x):=\eta(t,x,\vw_t(x))$. We have
    $${\mathcal L}:= -\int_{\R^{2d}}\eta_0(x,\vo)\,{\rm d}\mu_0 =  -\int_{\R^{d}}\eta_0(x,\vw_0(x))\,{\rm d}\vr_0=-\int_{\R^{d}}\widetilde{\eta}_0(x)\,{\rm d}\vr_0. $$
   
    Moreover
   \begin{equation}\label{eq:honk}
   \begin{split}
{\mathcal R}&:=\int_0^T\int_{\mathbb{R}^{2d}}\Big(\partial_t\eta_t(x,\omega) + (\omega-\nabla K*\vr_t(x))\cdot\nabla_x\eta_t(x,\omega)\Big) \,{\rm d}\mu_t \,{\rm d}t\\
&= \int_0^T\int_{\mathbb{R}^{d}}\Big(\partial_t\eta_t(x,\vw_t(x)) + (\vw_t(x)-\nabla K*\vr_t(x))\cdot\nabla_x\eta_t(x,\vw_t(x))\Big) \,{\rm d}\vr_t \,{\rm d}t
\end{split}
\end{equation}
We have
\begin{align*}
    \partial_t\widetilde{\eta}(t,x) &= \partial_t\eta(t,x,\vw_t(x)) + \nabla_\vw\eta(t,x,\vw_t(x))\partial_t\vw_t(x),\\
    \nabla_x\widetilde{\eta}(t,x) &= \nabla_x\eta(t,x,\vw_t(x)) + \nabla_\vw\eta(t,x,\vw_t(x))\nabla_x\vw_t(x).
\end{align*}
Thus, on the support of $\rho_t$, we have
\begin{align*}
    & \partial_t\widetilde{\eta}(t,x) + (\vw_t(x)-\nabla K*\vr_t(x))\cdot\nabla_x\widetilde{\eta}(t,x)\\
     &=\partial_t\eta(t,x,\vw_t(x)) + (\vw_t(x)-\nabla K*\vr_t(x))\cdot\nabla_x\eta_t(x,\vw_t(x))\\
     &\qquad+ \nabla_\vw\eta(t,x,\vw_t(x))\Big( \underbrace{\partial_t\vw_t(x) + \vu_t(x)\cdot\nabla\vw_t(x)}_{=0\ \mbox{ by } \eqref{AR_multi}} \Big).
\end{align*}

\noindent
Plugging the above into \eqref{eq:honk}, we infer that
$$ {\mathcal R}=\int_0^T\int_{\mathbb{R}^{d}}\Big(\partial_t\widetilde{\eta}_t(x) + (\vw_t(x)-\nabla K*\vr_t(x))\cdot\nabla_x\widetilde{\eta}_t(x)\Big) \,{\rm d}\vr_t \,{\rm d}t .$$
On the other hand, because $(\vr,\vw)$ is a weak solution to \eqref{AR_multi}, we have 
$$ \int_0^T\int_{\mathbb{R}^{d}}\Big(\partial_t\widetilde{\eta}_t(x) + (\vw_t(x)-\nabla K*\vr_t(x))\cdot\nabla_x\widetilde{\eta}_t(x)\Big) \,{\rm d}\vr_t \,{\rm d}t = -\int_{\R^{d}}\widetilde{\eta}_0(x)\,{\rm d} \vr_0 ,$$
which implies that ${\mathcal L} = {\mathcal R}$ and thus \eqref{weakeqF} holds true.  $\Box$

\begin{rem}
    The assumption that $\vw\in C^1_b([0,T]\times\R^d)$ can probably be weakened. One would need to use a mollification $\theta_\epsilon*\vw$ of $\vw$ in the composed test functions $\psi(x,\vw_t)$ and $\eta_t(x,\vw_t)$ in the proof of Proposition \ref{prop:monokin}. Then computing a commutator, say $\eta_t(x,\theta_\epsilon*\vw_t)-\eta_t(x,\vw_t)$ and passing to the limit as $\epsilon\to 0$, Proposition \ref{prop:monokin} can be recovered with weaker assumptions on the regularity of $\vw$. Similar strategy was successfully employed in \cite{PP22} in the proof of equivalence between the kinetic DNAR and Cucker-Smale models on the level of weak formulations.
\end{rem}

Let us conclude this section with a simple observation based on the results of Section \ref{sec:micromeso}, which allows us to connect the microscopic and the macroscopic DNAR systems. As a direct consequence of Proposition \ref{prop:monokin} and the results of Section \ref{sec:micromeso}, we obtain the following corollary.

\begin{cor}\label{corollary}
Suppose that $K$ is as in \eqref{defK} and in the $\lambda$-convex case it is $C^1$. Let 
$$0\leq\vr\in C([0,T];{\mathcal P}(\R^d)-\mbox{narrow}),\quad \vw\in C^1_b([0,T]\times\R^d)$$ be a weak solution to \eqref{AR_multi}. Then there exists a sequence of atomic solutions \eqref{E-empirical-measures} such that
\begin{equation*}
\sup_{t\in[0,T]}W_2(\mu_t^N, \vr(t,x)\delta(\vo-\vw))\xrightarrow{N\to\infty} 0 .
\end{equation*}
\end{cor}

\noindent
Let us take a closer look at the above corollary to better understand its implications.
For any bounded-Lipschitz function $\phi=\phi(x)$ we have
\begin{align*}
 \left|\int_{\R^{d}}\phi(x) d\vr^N_t - \int_{\R^{d}}\phi(x)\vr(t,x) dx\right| =   \left|\int_{\R^{2d}}\phi(x) d\mu^N_t - \int_{\R^{2d}}\phi(x) d\vr(t,x)\delta(\vo-\vw) \right|\\
    \leq d_{BL}(\mu_t^N, \vr(t,x)\delta(\vo-\vw))\leq W_2(\mu_t^N, \vr(t,x)\delta(\vo-\vw))\to 0,
\end{align*}
where, in analogy to previous notation $\vr^N_t = \int_{\R^d}\mu^N_t(x,d\vo)$.
Taking the supremum with respect to all $\phi$ with $|\phi|_\infty\leq 1$ and $[\phi]_{Lip}\leq 1$ we conclude that

\begin{equation}\label{mf1}
    \sup_{t\in[0,T]} d_{BL}(\vr^N_t,\vr(t,\cdot))\to 0. 
\end{equation}

\noindent
Similarly, denoting $\vr_t^N\vw_t^N:=\int_{\R^d}\omega\mu^N_t(x,\dd\omega)$, we can show that at each $t\in[0,T]$, we have
\begin{equation}\label{mf2}
    \int_{\R^d}\phi(x)\dd\vr_t^N\vw_t^N \to \int_{\R^d}\phi(x)\vw(t,x)\dd\vr_t(x),
\end{equation}

\noindent
for any bounded-Lipschitz test function $\phi$.
Indeed, by the boundedness of $\vw$ we may assume without loss of generality that for all $t\in[0,T]$ and all $N$ we have ${\rm supp}\, \mu^N_t\subset \R^d\times K$, where $K\subset \R^d$ is compact (the $\omega$-supports of the atomic solutions are uniformly bounded). Since the function $(x,\omega)\mapsto\phi(x)\omega$ is uniformly bounded and Lipschitz continuous in $\R^d\times K$ it is a good test function in the definition of the bounded-Lipschitz distance. Thus 
we have
\begin{align*}
\left|\int_{\R^{d}}\phi(x) \dd\vr^N_t\vw^N_t - \int_{\R^{d}}\phi(x)\vw(t,x)\dd\vr(t,x) \right| =    \left|\int_{\R^{2d}}\phi(x)\vo \dd\mu^N_t - \int_{\R^{2d}}\phi(x)\vo \dd\vr(t,x)\delta(\vo-\vw) \right|\\
   \leq d_{BL}(\mu_t^N, \vr(t,x)\delta(\vo-\vw)) \leq W_2(\mu_t^N, \vr(t,x)\delta(\vo-\vw))\to 0.
\end{align*}
Corollary \ref{corollary} and properties \eqref{mf1} and \eqref{mf2} require further investigation, but they provide the first step in the derivation of the micro-to-macroscopic limit for the DNAR system with kernels $K$ satisfying \eqref{defK}.

\section{Solutions to the DNAR system as a limit of particle dynamics}\label{micromacro}

The purpose of this section is to prove that, under stronger assumptions on the potential $K$, the classical solutions to the hydrodynamic DNAR system \eqref{AR_multi} might be obtained as a mean-field limit of the solutions to the microscopic formulation \eqref{ar}. We aim to  use the microscopic and macroscopic equivalence between the EAM and DNAR systems, see Proposition \ref{equiv} and subsequent comments.
Then we modify the results obtained by Carrillo-Choi in \cite{CC21} for the classical Euler-alignment system.

\subsection{The mean-field limit}
Another, currently better explored strategy follows the results of Carrillo and Choi \cite{CC21} concerning the direct micro-to-macroscopic mean-field limit for the pressureless Euler-alignment system. Below, we generalize them to the case of sufficiently smooth, matrix-valued communication weight $\boldsymbol{\Psi}\in W^{1,\infty}$. The matrix-valued equivalent of the mean-field limit result from \cite{CC21} (without the interaction potential) states the following.

\begin{theo}\label{CCth}
    Let $\{(x_i,v_i)\}_{i=1,\dots,N}$ solve the particle system \eqref{cs} and let
    $$(\vr,\vu)\in ( C([0,T];\mathcal{P}(\R^d)), L^\infty(0,T;W^{1,\infty})),\quad \varrho>0,$$ 
     be a classical solution to system \eqref{EA} subject to initial conditions $(\varrho,\vu)_{|_{t=0}} = (\varrho_0,\vu_0)$.


  \noindent
  Assume that the interaction matrix $\boldsymbol{\Psi}\colon\R^d\to\R^{d\times d}$ belongs to $W^{1,\infty}(\R^d)$, and for all $x\in\R^d$ satisfies
    \[ \boldsymbol{\Psi}(-x)=\boldsymbol{\Psi}(x),\qquad {\bf{v}}^T\boldsymbol{\Psi}(x){\bf{v}} \geq 0 \quad \forall_{{\bf{v}}\in \R^d}. \]

\noindent
    Then, the following estimate is satisfied:
    \begin{equation}\label{CC_estimate}
        \frac{1}{N}\sum_{i=1}^N |\vu(x_i(t),t)-v_i(t)|^2 + d_{BL}^2(\varrho^N,\varrho)
        \leq C\left(\frac{1}{N}\sum_{i=1}^N |\vu_0(x_i(0))-v_i(0)|^2 + d_{BL}^2(\varrho^N_0,\varrho)\right),
    \end{equation}
    where
    \[ \varrho^N 
    =\frac{1}{N}\sum_{i=1}^N\delta_{x_i(t)}. \]

    \noindent
    In consequence, assuming that the initial data satisfy 
    \begin{equation}\label{CC_initial} \lim_{N\to\infty}\left(\frac{1}{N}\sum_{i=1}^N |\vu_0(x_i(0))-v_i(0)|^2 + d_{BL}^2(\varrho^N_0,\varrho)\right) = 0, \end{equation}
    the following convergence holds:
   \eq{ \label{converg}     
        \frac{1}{N}\sum_{i=1}^Nv_i(t)\delta_{x_i(t)} &\rightharpoonup^* \varrho \vu \quad \text{in} \quad L^\infty(0,T;\mathcal{M}(\R^d)), \\
        \frac{1}{N}\sum_{i=1}^N v_i(t)\otimes v_i(t)\delta_{x_i(t)} &\rightharpoonup^* \varrho\vu\otimes\vu \quad \text{in} \quad L^\infty(0,T;\mathcal{M}(\R^d)), \\
        \frac{1}{N}\sum_{i=1}^N \delta_{x_i(t), v_i(t)} &\rightharpoonup^* \varrho\delta_{\vu} \quad \text{in} \quad L^\infty(0,T;\mathcal{M}(\R^d\times\R^d)).
    }
\end{theo}

{\emph{Sketch of the proof.}} Derivation of the estimate (\ref{CC_estimate}) is based on Gronwall's lemma. Differentiating the left-hand side in time (the dependence of time is omitted), we get
\[ \frac{1}{2}\frac{d}{dt}\frac{1}{N}\sum_{i=1}^N|\vu(x_i)-v_i|^2 = \frac{1}{N}\sum_{i=1}^N(\vu(x_i)-v_i)\cdot(\partial_t \vu(x_i)+\nabla_x \vu(x_i)v_i - \dot{v}_i). \]
Next, using momentum equations (\ref{cs}) and (\ref{EA}) we obtain 
\[\begin{aligned}
    &\frac{1}{2}\frac{d}{dt}\frac{1}{N}\sum_{i=1}^N|\vu(x_i)-v_i|^2\\
    =& -\frac{1}{N}\sum_{i=1}^N(\vu(x_i)-v_i)^T\nabla \vu(x_i)(\vu(x_i)-v_i) \\
    &+ \frac{1}{N}\sum_{i=1}^N(\vu(x_i)-v_i)^T\cdot\Big(\int_{\R^d}\boldsymbol{\Psi}(x_i-y)(\vu(y)-\vu(x_i))\varrho(t,y)\dy  -\frac{1}{N}\sum_{j\neq i}\boldsymbol{\Psi}(x_i-x_j)(v_j-v_i)\Big) \\
    =& I_1+I_2
\end{aligned}\]
For $I_1$ we get
\[ I_1 \leq \|\nabla \vu\|_{L^\infty([0,T]\times\R^d)}\frac{1}{N}\sum_{i=1}^N|\vu_i(x_i)-v_i|^2. \]
For the second part of $I_2$ we write
\[\begin{aligned} 
&-\frac{1}{N}\sum_{j\neq i}\boldsymbol{\Psi}(x_i-x_j)(v_j-v_i) \\
=& \frac{1}{N}\sum_{j\neq i}\boldsymbol{\Psi}(x_i-x_j)(\vu(x_j)-v_j-\vu(x_i)+v_i) - \frac{1}{N}\sum_{j\neq i}\boldsymbol{\Psi}(x_i-x_j)(\vu(x_j)-\vu(x_i)) \\
=& -\frac{1}{N}\sum_{j\neq i}\boldsymbol{\Psi}(x_i-x_j)(\vu(x_j)-v_j-\vu(x_i)+v_i) - \int_{\R^d}\boldsymbol{\Psi}(x_i-y)(\vu(y)-\vu(x_i))\;\mathrm{d}\varrho^N(y).
\end{aligned}\]
Using the above, we express ${I_2}$ as
\[\begin{aligned} {I_2} =& \frac{1}{N}\sum_{i=1}^N(\vu(x_i)-v_i)^T\int_{\R^d}\boldsymbol{\Psi}(x_i-y)(\vu(y)-\vu(x_i))\;\mathrm{d}(\varrho-\varrho^N)(y) \\
&+ \frac{1}{N^2}\sum_{i,j}(\vu(x_i)-v_i)^T\Psi(x_i-x_j)(\vu(x_j)-v_j-\vu(x_i)+v_i).
\end{aligned}\]
However, by the symmetry of $\boldsymbol{\Psi}$ the second part of $I_2$ equals
\[ -\frac{1}{N^2}\sum_{i,j}(\vu(x_i)-v_i-\vu(x_j)+v_j)^T\boldsymbol{\Psi}(x_i-x_j)(\vu(x_i)-v_i-\vu(x_j)+v_j), \]
which is negative by the fact that $\boldsymbol{\Psi}$ is positive semi-definite, and thus, we arrive at
\[ {I_2} \leq \frac{1}{N}\sum_{i=1}^N(\vu(x_i)-v_i)\int_{\R^d}\boldsymbol{\Psi} (x_i-y)(\vu(y)-\vu(x_i))\;\mathrm{d}(\varrho-\varrho^N)(y). \]
Finally, applying H\"older's inequality we obtain
\[ I_2 \leq \frac{1}{N}\sum_{i=1}^N|\vu(x_i)-v_i|^2 + \||\boldsymbol{\Psi}| |\vu|\|_{W^{1,\infty}}^2d_{BL}^2(\varrho,\varrho^N). \]
From Proposition 2.2 from \cite{CC21}, it follows that
\[ d_{BL}(\varrho,\varrho^N) \leq d_{BL}(\varrho_0,\varrho_0^N) + C\sqrt{\frac{1}{N}\sum_{i=1}^N|\vu(x_i)-v_i|^2}. \]
Note that this estimate is derived only using the continuity equation for $\varrho$ and the definition of $\varrho^N$, and thus remains the same for the matrix-valued communication weight. In conclusion 
\[ \frac{\mathrm{d}}{\mathrm{d}t}\frac{1}{N}\sum_{i=1}^N|\vu(x_i)-v_i|^2 \leq C\left(\frac{1}{N}\sum_{i=1}^N|\vu( x_i)-v_i|^2 + d^2_{BL}(\varrho_0,\varrho_0^N)\right) \]
and the estimate (\ref{CC_estimate}) follows from Gronwall's lemma.

From the estimate (\ref{CC_estimate}), one can easily deduce the desired convergence. Indeed, for any $\boldsymbol{\varphi}\in \lr{W^{1,\infty}(\R^d)}^d$, we have 
\[\begin{aligned} 
\left|\int_{\R^d}\boldsymbol{\varphi}(x)\cdot(\varrho \vu)(x)\;\mathrm{d}x-\frac{1}{N}\sum_{i=1}^N v_i\cdot\boldsymbol{\varphi}(x_i)\right| 
=& \left|\frac{1}{N}\sum_{i=1}^N\boldsymbol{\varphi}(x_i)\cdot(\vu(x_i)- v_i) + \int_{\R^d}\boldsymbol{\varphi}(x)\cdot\vu(x)\;\mathrm{d}(\varrho-\varrho^N)(x)\right| \\
\leq & \|\boldsymbol{\varphi}\|_{L^\infty}\left(\frac{1}{N}\sum_{i=1}^N|\vu(x_i)- v_i|^2\right)^{1/2} + \|\boldsymbol{\varphi}\cdot\vu\|_{W^{1,\infty}}d_{BL}(\varrho,\varrho^N).
\end{aligned}\]
Similarly, for any $\boldsymbol{\varphi}\in \lr{W^{1,\infty}(\R^d)}^{d\times d}$, we have
\[\begin{aligned} 
&\int_{\R^d}\boldsymbol{\varphi}:(\varrho\vu\otimes\vu)\;\mathrm{d}x - \frac{1}{N}\sum_{i=1}^N\boldsymbol{\varphi}(x_i): v_i\otimes v_i \\
=& -\frac{1}{N}\sum_{i=1}^N\boldsymbol{\varphi}(x_i):(\vu(x_i)- v_i)\otimes(\vu(x_i)- v_i) 
+ \int_{\R^d}\boldsymbol{\varphi}:(\vu\otimes \vu)\;\mathrm{d}(\varrho-\varrho^N)(x) \\
&+ \frac{1}{N}\sum_{i=1}^N\boldsymbol{\varphi}(x_i):\Big(\vu(x_i)\otimes(\vu(x_i)- v_i) +(\vu(x_i)- v_i)\otimes\vu(x_i)\Big),
\end{aligned}\]
and so
\[\begin{aligned}
    \left|\int_{\R^d}\boldsymbol{\varphi}:(\varrho\vu\otimes\vu)\;\mathrm{d}x - \frac{1}{N}\sum_{i=1}^N\boldsymbol{\varphi}(x_i): v_i\otimes v_i\right| \leq & \|\boldsymbol{\varphi}\|_{L^\infty}\frac{1}{N}\sum_{i=1}^N|\vu(x_i)- v_i|^2 + \||\boldsymbol{\varphi}||\vu|^2\|_{W^{1,\infty}}d_{BL}(\varrho,\varrho^N) \\
    &+ 2\||\boldsymbol{\varphi}||\vu|\|_{L^\infty}\left(\frac{1}{N}\sum_{i=1}^N|\vu(x_i)- v_i|^2\right)^{1/2}
\end{aligned}\]
In conclusion, assuming that the initial data satisfy (\ref{CC_initial}), we have convergences \eqref{converg} with respect to the bounded-Lipshitz distance. 

On the other hand, for the empirical measure $\mu_t^N(x,v)=\frac{1}{N}\sum_{i=1}^N \delta_{(x_i(t), v_i(t))} $, and a test function ${\varphi}\in W^{1,\infty}(\R^{2d})$, we have
\[\begin{aligned}
    &\left|\int_{\R^d\times\R^d}\varphi(x, v)\;\mathrm{d}(\varrho\delta_{\vu}-\mu^N)(x, v)\right| \\=& \left|\int_{\R^d}\varphi(x,\vu(x))\;\mathrm{d}(\varrho-\varrho^N)(x) + \frac{1}{N}\sum_{i=1}^N(\varphi(x_i,\vu(x_i))-\varphi(x_i, v_i))\right| \\
    \leq & \|\varphi\|_{W^{1,\infty}}d_{BL}(\varrho,\varrho^N) 
    + (\|\varphi\|_{L^\infty}+\|\varphi\|_{W^{1,\infty}}\|\vu\|_{W^{1,\infty}}\left(\frac{1}{N}\sum_{i=1}^N|\vu(x_i)- v_i|^2\right)^{1/2},
\end{aligned}\]
and so 
\[ d_{BL}(\varrho\delta_{\vu},\mu^N)\to 0. \]

\subsection{Existence of strong solutions}
In this section we state the result for the existence of local in time classical solutions to the hydrodynamic DNAR system \eqref{AR_multi}. We again  employ the observation from the introduction, that for  classical solutions this system is equivalent with the EAM system \eqref{EA}. The following statements generalize the results of \cite{CC21}, to the case of matrix-valued communication weight.

\begin{defi}\label{df:reg}
Let 
$s>\frac{d}{2}+1$. For given $T\in (0,\infty)$, the pair $\vr,\vu$ is a strong solution  of \eqref{EA} with \eqref{khess} on the interval $[0,T]$ if the following conditions are satisfied:
\begin{enumerate}
\item $\vr\in C([0,T];H^s(\R^d))$, $\vu\in C([0,T];Lip (\R^d)\cap L^2_{{\rm{loc}}}(\R^d))$ and $\Grad^2 \vu\in C([0,T];H^{s-1}(\R^d))$,
\item $(\vr,\vu)$ satisfy equations \eqref{EA} in the sense of distributions.
\end{enumerate}
\end{defi}

The analogue of the Theorem 4.1 from \cite{CC21} in the case of matrix-valued communication weight reads as follows.

\begin{theo}\label{Th:local}
Let $s>\frac{d}{2}+1$ and $R>0$. Suppose that the communication weight $\boldsymbol{\Psi}$ satisfy
\eqh{
\boldsymbol{\Psi}\in C^1_c(\R^d)\quad and \quad supp(\boldsymbol{\Psi})\subseteq B(0,R),
}
where $B(0,R)\subset\R^d$ denotes a ball of radius $R$ centered at origin. For any $N<M$, there exists a positive constant $T^*$ depending only on $R,N,$ and $M$ such that if $\vr_0>0$ and
\eqh{
\|\vr_0\|_{H^s}+\|\vu_0\|_{L^2(B(0,R))}+\|\Grad \vu_0\|_{L^\infty}+\|\Grad^2 \vu_0\|_{H^{s-1}}<N,
}
then the Cauchy problem \eqref{EA} has a unique strong solution $(\vr,\vu)$, in the sense of Definition \eqref{df:reg}, satisfying
\eqh{
\sup_{0\leq t\leq  T^*}\lr{\|\vr(\cdot,t)\|_{H^s}+\|\vu(\cdot,t)\|_{L^2(B(0,R))}+\|\Grad \vu(\cdot,t)\|_{L^\infty}+\|\Grad^2 \vu(\cdot,t)\|_{H^{s-1}}}\leq M.
}
\end{theo}

Instead of essentially repeating the proof from \cite{CC21}, let us conclude with the following remarks.
\begin{rmk}
 The compactness of the support is necessary for this construction to work. This assumption could be relaxed on the bounded domain. 
\end{rmk}

\begin{rmk}
In comparison to \cite{CC21} the proof of Theorem \eqref{Th:local} does not require additional estimate of local-in-space norms. This is due to lack of attraction-repulsion type of interactions in the system.
\end{rmk}

\appendix
\section{Disintegration Theorem}

For the reader's convenience we recall the disintegration theorem as stated in \cite[Theorem 5.3.1]{AGS}.

\begin{theo}[\textbf{Disintegration theorem}]\label{T:disint}
For $d_1, d_2\in\mathbb{N}$ denote projection onto the second componenent as $\pi_2:\R^{d_1}\times \R^{d_2} \longrightarrow \R^{d_2}$. For $\mu \in \mathcal{P} \left( \R^{d_1} \times \R^{d_2}\right)$, define its projection onto second factor as $\nu :=(\pi_2)_{\#} \mu \in \mathcal{P}(\R^{d_2})$. Then there exists a family of probabilistic measures $\left\{ \mu_{x_2}\right\}_{x_2 \in \R^{d_2}}\subset \mathcal{P}(\R^{d_1})$, defined uniquely $\nu$ almost everywhere, such that 
\begin{enumerate}
    \item The map
    $$
    \R^{d_2} \ni x_2 \longmapsto \mu_{x_2}(B)
    $$
    is Borel-measurable for each Borel set $B \subset \R^{d_1}$;
    \item The following formula 
    $$
    \int_{\R^{d_1}\times\R^{d_2}} \phi(x_1, x_2) \dd \mu(x_1,x_2) = \int_{\R^{d_2}}\left( \int_{\R^{d_1}} \phi(x_1, x_2) \dd \mu_{x_2}(x_1) \right) \dd \nu (x_2)
    $$
    holds for every Borel-measurable map $\phi:\R^{d_1} \times \R^{d_2} \longmapsto [0, \infty)$.
\end{enumerate}
\end{theo}

\noindent{\bf{Acknowledgements}}~ 
The research of N.C. was supported by the EPSRC Early Career Fellowship EP/V000586/1, and by the “Excellence Initiative Research University (IDUB)” program at the University of Warsaw. J.P.'s work has been supported by the Polish National Science Centre’s Grant No.
2018/31/D/ST1/02313 (SONATA) and partially by University of Warsaw
program IDUB Nowe Idee 3A. The work of M. S. was supported by the National Science Centre grant no. 2022/45/N/ST1/03900 (Preludium). The work of N.C. and E.Z. was  supported by the EPSRC Early Career Fellowship no. EP/V000586/1.


\begin{thebibliography}{99}

\bibitem{ABDM}
 P. Aceves-S\'{a}nchez, R. and Bailo, P. Degond,  and Z. Mercier.
\newblock Pedestrian models with congestion effects. 
\newblock{\em Math. Models Methods Appl. Sci.}, 1--41, doi: 10.1142/S0218202524400050, 2024.
\bibitem{AGS}
L. Ambrosio, N. Gigli, and G. Savaré.  
\newblock {\em{Gradient flows in metric spaces and in the space of probability measures}.}
\newblock Birkhäuser, Basel, 2008.
	
    	\bibitem{AR}
    	A. Aw and M. Rascle. 
    	\newblock Resurrection of second order models of traffic flow.
    	\newblock {\em SIAM J. Appl. Math.}, 60:916--938, 2000.

 \bibitem{CC21}
J. A. Carrillo and Y.-P. Choi.
\newblock Mean-field limits: from particle descriptions to macroscopic equations.
{\it Arch. Ration. Mech. Anal.}, {\bf 241}, 1529--1573, 2021. 

\bibitem{CGZ} 
N. Chaudhuri, P. Gwiazda and E. Zatorska.
\newblock Analysis of the generalised Aw-Rascle model. 
\newblock {\em Comm. Partial Differential Equations} {\bf 48}, no. 3, 440–477, 2023.

	\bibitem{CMPZ}  N. Chaudhuri, M. Ali Mehmood, C. Perrin and E. Zatorska.
	\newblock  Duality solutions to the hard-congestion model for the dissipative Aw-Rascle system. 
	\newblock  {\em Comm. Partial Differential Equations.} https://doi.org/10.1080/03605302.2024.2380696, 2024.
\bibitem{CNPZ}
N. Chaudhuri, L. Navoret, C. Perrin and E. Zatorska.
\newblock Hard congestion limit of the dissipative Aw-Rascle system.
	\newblock {\em Nonlinearity.} {\bf 37}, 04, 2024.
\bibitem{CZ}
Y-P. Choi, and X. Zhang,
\newblock One dimensional singular Cucker-Smale model: uniform-in-time mean-field limit and contractivity.
{\it J. Differential Equations.}, {\bf 287}, 428–459, 2021.

     \bibitem{CDS20}
P. Constantin, T. D. Drivas, and R. Shvydkoy. 
\newblock Entropy Hierarchies for equations of compressible fluids and self-organized dynamics. 
{\em SIAM J. Math. Anal.}, {\bf{52}}, 3073--3092, 2020.

     \bibitem{kis}
T. Do, A. Kiselev, L. Ryzhik, and C. Tan.
\newblock Global regularity for the fractional {E}uler alignment system.
\newblock {\em Arch. Ration. Mech. Anal.}, 228(1):1--37, 2018.

\bibitem{FP} { M. Fabisiak and J. Peszek, Inevitable monokineticity of strongly singular alignment, {\it Mathematische Annalen},  https://doi.org/10.1007/s00208-023-02776-7, 2023.}

\bibitem{FK}
A. Figalli and M.-J. Kang.
\newblock 
A rigorous derivation from the kinetic Cucker–Smale model to the pressureless Euler system with nonlocal alignment.
\newblock{\emph{Analysis \& PDE}}, 12(3), 843–-866, 2019.

\bibitem{HKPZ}
S.-Y. Ha, J. Kim, J. Park, and X. Zhang.
\newblock Complete cluster predictability of the Cucker-Smale flocking model on the real line.
\newblock {\emph{Arch. Ration. Mech. Anal.}}, 231, no. 1, 319–--365, 2019.

\bibitem{HPZ}
S.-Y. Ha, J. Park, and X. Zhang. 
\newblock A first-order reduction of the Cucker-Smale model on the real line and
its clustering dynamics. 
\newblock {\emph{Commun. Math. Sci.}}, 16, 1907–-1931,  2018.

\bibitem{HT}
{S.-Y. Ha and E. Tadmor, From particle to kinetic and hydrodynamic descriptions of flocking, {\it Kinet. Relat. Models}, {\bf 1} , no.~3, 415--435, 2008.
}

\bibitem{KMT}{T.~K. Karper, A. Mellet and K. Trivisa, Hydrodynamic limit of the kinetic Cucker-Smale flocking model, {\it Math. Models Methods Appl. Sci.} {\bf 25}, no.~1, 131--163, 2015.
}

\bibitem{Kim}
J. Kim. 
\newblock A Cucker-Smale flocking model with the hessian communication weight and its first-order reduction. 
\newblock {\em J. Nonlinear Sci.} 32(20),  2022.

\bibitem{Kim21}
J. Kim. 
\newblock First-order reduction and emergent behavior of the one-dimensional kinetic Cucker-Smale equation,
\newblock {\em J. Differ. Equations}, 302, 496-–532, 2021.

\bibitem{LS}
{T. M. Leslie  and R. Shvydkoy}.
\newblock {On the structure of limiting flocks in hydrodynamic Euler Alignment models},
\newblock{\em Math. Models Methods Appl. Sci.}, 29(13), {2419--2431}, {2019}.

\bibitem{M}
M. Ali Mehmood.
\newblock Hard congestion limit of the dissipative Aw-Rascle system with a polynomial offset function.
\newblock {\em J. Math. Anal. Appl.} 533, no. 1, 2024.

\bibitem{PP22}
J. Peszek and D. Poyato. 
\newblock Measure solutions to a kinetic Cucker-Smale model with singular and matrix-valued communication.
\newblock arXiv:2207.14764.

\bibitem{PP23}
J. Peszek and D. Poyato. 
\newblock Heterogeneous gradient flows in the topology of fibered optimal transport. 
\newblock{\em Calc. Var.} 62:258, 2023. 

\bibitem{ST21}
R. Shu and E. Tadmor. 
\newblock Anticipation Breeds Alignment. 
\newblock {\em Arch. Rat. Mech. Anal.}, 240:203–241, 2021.

\bibitem{tad1}
R. Shvydkoy and E. Tadmor.
\newblock Eulerian dynamics with a commutator forcing.
\newblock {\em Transactions of Mathematics and Its Applications}, 1(1):tnx001,
  2017.

\bibitem{tad2}
R.~Shvydkoy and E.~Tadmor.
\newblock Eulerian dynamics with a commutator forcing {II}: flocking.
\newblock {\em Disc. and Cont. Dyn. Sys.}, 37(11):5503--5520, 2017.



\bibitem{tad4}
R. Shvydkoy and E. Tadmor.
\newblock Eulerian dynamics with a commutator forcing {III}. fractional
  diffusion of order $0<\alpha<1$.
\newblock {\em Physica D: Nonlinear Phenomena}, 376-377:131 -- 137, 2018.
\newblock Special Issue: Nonlinear Partial Differential Equations in
  Mathematical Fluid Dynamics.

     \bibitem{Zhang} M. Zhang.
\newblock A non-equilibrium traffic model devoid of gas-like behavior.
\newblock {\em Transportation Res. B} 36, 275--290, 2002.







\end{thebibliography}
\end{document}